\RequirePackage{etoolbox}
\csdef{input@path}{%
 {sty/}
 {img/}
}%
\csgdef{bibdir}{bib/}

\documentclass{article}

\newcommand{\nb}{\mbox{NB}}
\newcommand{\E}{\mbox{E}}
\newcommand{\x}{\bm X}

\usepackage{amsthm}
\usepackage{amsmath}
\usepackage{graphicx}
\usepackage{bm}
\usepackage{amsfonts}
\usepackage{caption}
\usepackage{tikz}
\usepackage{subcaption}
\usepackage{url}

\author{Anna Heath, Ioanna Manaolopoulou and Gianluca Baio}

\begin{document}


\title{Efficient Monte Carlo Estimation of the Expected Value of Sample Information using Moment Matching}


\maketitle

\begin{abstract}
The Expected Value of Sample Information (EVSI) is used to calculate the economic value of a new research strategy. While this value would be important to both researchers and funders, there are very few practical applications of the EVSI. In the main, this is due to computational difficulties associated with calculating the EVSI in practical health economic models using nested simulations. We present an approximation method for the EVSI that is based on estimating the distribution of the posterior mean of the incremental net benefit across all the possible future samples, known as the distribution of the preposterior mean. Specifically, we suggest that this distribution is estimated using moment matching coupled with simulations that are available for probabilistic sensitivity analysis, which is typically mandatory in health economic evaluation. We demonstrate that this method is successful using an example that has previously been applied to other EVSI approximation methods. We then conclude by discussing how our method fits in with other recent additions to the literature that detail approximation methods for the EVSI. 
\end{abstract}

Value of Information (VoI) analysis \cite{Howard:1966} has frequently been touted as an extension to probabilistic sensitivity analysis (PSA) \cite{FelliHazen:1998,FelliHazen:1999,Claxton:1999,Claxtonetal:2001,Adesetal:2004,
BrennanKharroubi:2005,Briggsetal:2006,Fenwicketal:2006}. Generally speaking, VoI analyses quantify the potential economic benefit of performing new research targeting uncertainty in an underlying economic model. In the past, these measures have focused on estimating the economic benefit of resolving \emph{all} the uncertainty in the model. While this is simple to calculate computationally, it has limited practical value, especially in models where the benefit is high, as model uncertainty can almost never be fully resolved.

Arguably, therefore, the VoI measure that is most useful in practice is known as the Expected Value of \emph{Sample} Information (EVSI) \cite{Schlaifer:1959,RaiffaSchlaifer:1961}, which measures the potential economic benefit of some specific future research. This value derives directly from the possibility that the new information gleaned from the research could help a payer avoid wasting resources on an inefficient treatment which would, however, have been deemed cost-effective given the current knowledge about the underlying model inputs.

Despite the obvious benefits of calculating the EVSI, especially in terms of efficiently using research and development budgets, practical applications of EVSI calculations are scarce \cite{Steutenetal:2013}. Initially, this was because the published approaches to calculating the EVSI were based on computationally intensive nested Monte Carlo (MC) procedures \cite{Brennanetal:2007}. Some older methods avoid the need for nested simulation by assuming certain model conditions such as independence in the model parameters and linearity \cite{Adesetal:2004}, which reduced the computation time but limited the applicability of these methods.

More recently, EVSI research has been active and this has led to the proposal of several alternative calculation methods. Some of these more recent methods calculate the EVSI for specific designs, such as cluster randomised clinical trials \cite{Weltonetal:2014,BrennanKharroubi:2007b,Adesetal:2004}, while others use approximations to avoid using nested simulations \cite{BrennanKharroubi:2007,Kharroubietal:2011,Jalaletal:2015}. The most recent approaches have been developed to avoid rerunning the economic model, which means that the EVSI can be calculated using the samples that have already been obtained to perform the, often mandatory, probabilistic sensitivity analysis (PSA) \cite{Menzies:2016,Strongetal:2015}. These methods may be suitable in certain situations, especially in settings where each model run requires a large amount of computational time. However, they either require the specification of sufficient statistics or a large number of matrix calculations.

In this paper, we present a new computation method for the EVSI that can be used irrespective of the underlying structure in the health economic model. The method is based on determining the properties of a specific function of the net benefit. We suggest an approximation based on these properties and then demonstrate how this can be used to calculate the EVSI. Specifically, the approximation reuses PSA samples from the underlying model alongside a small number of estimates of the variance of the incremental net benefit obtained by rerunning the model. Therefore, this method relies on around 30 nested simulations compared to at least 600 \cite{Oakleyetal:2010} for a standard EVSI analysis. 

In models with a relatively short run time, our method has a similar computational cost to more recent alternatives. It is also simple to implement for virtually all model structures and realistic trial designs, possibly accounting for missingness or problems with follow-up, as it relies simply on Bayesian updating and variance calculations. Most importantly, the required Bayesian updating is in the same form as the analysis that would be required once the data are collected, which aligns the model used to plan the study with the one used for the actual data analysis. This is not always the case with traditional methods, when sample size calculations are based on simplistic approximations and analytic formul\ae.

The presentation of our method begins in \S\ref{Notation} with the introduction of notation and other key concepts. The estimation method is then presented in \S\ref{pre-post} and it is implemented for an example from Ades \textit{et al.} \cite{Adesetal:2004} in \S\ref{toy}. Finally, we conclude with a discussion of the computational burden of this method and how this compares to other methods. This allows for the recommendation of when this new tool may be more suitable than the alternatives. 

\section{Notation and Concepts}\label{Notation}
The EVSI is calculated as the difference between the value of the decision made under the current level of uncertainty and the expected value of the decision made with the additional information contained in the future trial. To present both the EVSI and our method more formally, we begin by introducing some notation and key concepts. First, we assume that the health economic model is defined using a set of parameters, denoted $\bm\theta$. To perform Probabilistic Sensitivity Analysis (PSA), the current level of uncertainty in these parameters, based on  literature reviews, clinical trials or meta analyses, is defined using a probability distribution $p(\bm\theta)$ \cite{BaioDawid:2011}.

To calculate the EVSI, it is assumed that $T$ treatments are under consideration. To compare these treatments, we define the \emph{net benefit} of each treatment $t=1,\dots,T$ based on these model parameters, which we indicate as $\nb_t^{\bm\theta}$. Strictly speaking, we assume that $\nb_t^{\bm\theta}$ does not reflect individual level uncertainty so if all the $\bm\theta$ values were known with certainty then $\nb_t^{\bm\theta}$ would be a fixed number. The \emph{value} of the best decision under current information is the treatment with the maximum net benefit value \[
\max_t \E_{\bm\theta}\left[\nb_t^{\bm\theta}\right].
\]

The EVSI is concerned with a trial in which we gain ``more information'' by gathering a new dataset $\bm X$. This new dataset is what we would collect in the future trial, for example the number of patients responding to a treatment or the number of false positive test results. These data are useful as they provide us with more information about the underlying model parameters which, in turn, influence the decision. Specifically, if the future trial had already been undertaken and the data observed as $\bm x$, the optimal decision conditional on that observed sample would be \[\max_t \mbox{E}_{\bm\theta\mid\bm X=\bm x}\left[\mbox{NB}_t^{\bm\theta}\right],\] where the expectation of the net benefit is taken over the distribution of the parameters \emph{conditional} on the information in the sample $\bm x$.

However, as new research has not been carried out yet and the data $\bm X$ have not been observed, the EVSI is calculated by averaging over all possible future datasets to give the average value of the decision made with the additional information contained in the sample; \begin{equation}\mbox{EVSI} = \mbox{E}_{\bm X}\left[\max_t \mbox{E}_{\bm\theta\mid\bm X}\left[\nb_t^{\bm\theta}\right]\right] - \max_t\mbox{E}_{\bm\theta}\left[\mbox{NB}_t^{\bm\theta}\right].\label{EVSI}\end{equation} In general, the distribution of all the possible future samples is defined using the same distribution that would be used to model the data if it had been observed \cite{Strongetal:2015,Menzies:2016}. Formally, this means that the distribution of $\bm X$ is defined through its relationship with $\bm\theta$, this allows us to sample values from the distribution of $p(\x)$ as we will see in \S\ref{pre-post var}. 

The analysis required for the EVSI is exactly the analysis that would be required if the statistical modelling underpinning the health economic model was inferred using Bayesian methods. In this setting, the PSA distributions for the parameters would be called the \emph{prior} for $\bm\theta$ and this \emph{prior} is combined with the data $\x$ to update the information about the model parameters and determine a \emph{posterior} distribution for the parameters $\bm\theta$ \emph{conditional} on the data $\bm X$. In fact, in Bayesian analysis, the inner expectation in the first term of equation (\ref{EVSI}) is known as the preposterior mean, i.e.~the posterior mean \emph{before} the data have been collected. The preposterior mean, which we will denote $\mu^{\x}_t$, is the key element in the EVSI calculation and it is also the term that makes the EVSI computationally expensive, as it has traditionally been estimated by simulation. In the most general setting, this requires posterior updating, usually with Markov Chain Monte Carlo (MCMC) methods, which are very computationally intensive. 

Therefore, all the approximation methods \cite{Adesetal:2004,BrennanKharroubi:2007,Jalaletal:2015,Strongetal:2015,Menzies:2016} focus on estimating this preposterior mean without using full MCMC sampling. Our approximation reuses the information contained in $\nb_t^{\bm\theta}$ to estimate the preposterior mean or more specifically the \emph{distribution} of the preposterior mean. The concept of ``distribution of a mean'' is rather counterintuitive in a standard statistical analysis. However, it makes sense in this context, because we are interested in the posterior \textit{before the data have been collected}. This implies that there is uncertainty as to which of all the possible future samples will occur, were we to conduct the study. Each possible sample would yield a different posterior mean and so the sampling distribution for the data induces a distribution over the posterior means. To make this more clear, we use a simple example, which is also used to highlight how $\nb_t^{\bm\theta}$ can be used to approximate this distribution.

\subsection{The distribution of the preposterior mean -- an example}\label{preposterior mean}\label{Example}
Suppose that a new drug is available and is associated with a probability $\theta$ of curing a particular disease. As this drug is new, we assume that there is very limited evidence on its effectiveness. This could be expressed, in a very simplistic way, by assuming that all values for $\theta$, between 0 and 1, are equally likely. Mathematically, this is equivalent to modelling $\theta\sim\mbox{Uniform}(0,1)$. Of course, in practical settings, it is likely that information about the effectiveness of the treatment is available and this could be easily included by using an alternative distribution for $p(\theta)$.

In this simple model, the treatment is effective if the disease has been cured and the drug costs are known to be equal to $c$. As we are only interested in \emph{population level} effectiveness, the effectiveness of the treatment is the probability of a cure $\theta$. We finally assume that the decision maker is prepared to pay $k$ monetary units for each person that is cured of the disease.

To simplify this example further, the new treatment is being compared with the current standard of care, which is to leave the disease untreated. This option has no cost but also has no effectiveness as this (non-life-threatening) disease does not improve without drug intervention. Thus, the two net benefit values are \[\mbox{NB}_1^{\theta} = 0 \quad \mbox{ and } \quad \mbox{NB}_2^{\theta} = k\theta-c.\]

The trial involves giving $N$ people the drug and observing how many are cured. This means that the sampling variability can be characterised using a Binomial distribution $X\mid\theta \sim \mbox{Binomial}(N,\theta),$  with $\theta$ the probability of being cured. The distribution $p(X\mid\theta)$ and the \emph{prior} for $\theta$ can then be combined to give the predictive distribution of the future samples $X$ as \begin{align*}p(X) &= \int_0^1 p(X\mid\theta)p(\theta)\, d\theta \\&= \int_0^1 {N \choose X} \theta^X (1-\theta)^{N-X} \times 1 \,d\theta \\&=\frac{X!(N-X)!N!}{(N+1)!X!(N-X)!}=\frac{1}{N+1}.\end{align*} This calculation implies that, if we begin by assuming that all values for $\theta$ are equally likely, then all values of the future sample are equally likely. This is because the possible values for $X$, the number of people cured, are $0,1,\dots,N$ and each of them has probability $\frac{1}{N+1}$.

Once the distribution for the possible data values is known, the distribution of the preposterior mean is found by calculating the posterior mean for both the net benefit functions. As $\nb_1^\theta$ does not depend on $\theta$, the posterior mean for NB$_1^{\theta}$ is $\mu^{X}_1=0$. Therefore, the preposterior mean for $\nb_1^\theta$ does not have a distribution in this example and is just equal to 0. However, the posterior mean for NB$_2^{\theta}$ does depend on the future sample $X$: \begin{align*}
\mu^{ X}_2= \mbox{E}_{\theta\mid X}\left[\mbox{NB}_2^{\theta}\right] &= \int_0^1 (k\theta-c)\ p(\theta\mid X)\, d\theta\\&=k\left(\int_0^1 \theta\, p(\theta\mid X)\,d\theta\right)-c \\
&= k\frac{1+ X}{2+N}-c ,
\end{align*}
as the second integral is exactly the mean of $\theta\mid X$ and it can be shown that $\theta\mid X \sim \mbox{Beta}(1+ X, 1+N)$. Therefore, the distribution of $\mu_2^X$ is directly related to the distribution of $X$, which assumes that all values for $X$ are equally likely and thus all values of $\mu_2^X$ are equally likely. However, this belief for $X$ is a direct result of our initial beliefs that all values of $\theta$ are equally likely. In this sense, the distribution of the preposterior mean is strongly linked to our initial beliefs about $\theta$ which are encoded in our prior distribution. The supplementary material presents two alternative models in order to demonstrate that the prior distribution for $\bm\theta$ is similar to the distribution of the preposterior mean, in many settings.

\section{Estimating the distribution of the preposterior mean}\label{pre-post}

Historically, EVSI calculations have been carried out by estimating the distribution of the preposterior mean by Monte Carlo simulation \cite{Prattetal:1995,Brennanetal:2007,Adesetal:2004} or by approximating the relationship between the samples $\bm X$ and the posterior mean \cite{Strongetal:2015,Adesetal:2004}. However, we have just demonstrated that the distribution of the preposterior mean for the net benefit, which we will indicate with $p(\mu_{t}^{\bm X})$, is related to the prior for the model parameters. More importantly, in health economic evaluations, PSA samples from the \emph{prior} are already available when the EVSI is being calculated \cite{BaioDawid:2011,Baio:2012,Andronisetal:2009}. Therefore, these PSA samples can be used to save computational time when estimating $p(\mu_{t}^{\bm X})$\footnote{Menzies \cite{Menzies:2016} also uses the similarity between the PSA samples for the net benefit and the distribution of the preposterior mean as a basis for his method but the proposed transformations differ significantly.}. Specifically, our method estimates $p(\mu_{t}^{\bm X})$ using the PSA samples and additional knowledge about the mean and variance of the preposterior mean, which are both estimated by simulation.

\subsection{Expectation and Variance for the preposterior mean}\label{Moment}
To estimate the mean and variance of $\mu_{t}^{\bm X}$, we use formul\ae\ for conditional iterated expectation \cite{Weiss:2006}. Firstly, this implies that the mean of $\mu_t^{\bm X}$ is given by \[\mbox{E}_{\bm X}\left[\mu^{\bm X}_t\right] = \mbox{E}_{\bm X}\left[ \mbox{E}_{\bm\theta\mid\bm X}\left[\mbox{NB}_t^{\bm\theta}\right] \right] = \mbox{E}_{\bm\theta}\left[\mbox{NB}_t^{\bm\theta}\right],\] i.e.~the average of the preposterior mean is equal to the mean of the $\nb_t^{\bm\theta}$, which can be obtained using the PSA samples for the net benefit. This implies that performing the EVSI analysis, i.e.~considering the possible values of a future data set before obtaining data, cannot give any additional information or change our current decision, as on average the optimal decision would remain the same. To actually change the optimal decision the research must be carried out.

Secondly, the variance of the preposterior mean is given by \[\mbox{Var}_{\bm X}\left[\mu^{\bm X}_t\right] = \mbox{Var}_{\bm X}\left[\mbox{E}_{\bm\theta\mid\bm X}\left[\mbox{NB}_t^{\bm\theta}\right]\right] = \mbox{Var}_{\bm\theta}\left[\mbox{NB}_t^{\bm\theta}\right] - \mbox{E}_{\bm X}\left[\mbox{Var}_{\bm\theta\mid\bm X}\left[\mbox{NB}_t^{\bm\theta}\right]\right].\] In other words, the variance of the preposterior mean is equal to the variance of the net benefit minus the average posterior variance, over all possible future samples. 

Therefore, to calculate the mean and variance of the preposterior mean by simulation, we only need to estimate the average posterior variance over all possible samples $\bm X$, as long as the PSA samples have already been obtained. Most importantly, \S\ref{pre-post var} demonstrates that a suitable estimate of the expected posterior variance can be obtained using only a small number of simulated future samples. This significantly reduces the number of samples required to calculate the EVSI compared to standard simulation based methods.

\subsection{Moment Matching}\label{Moment}
Assuming that the mean and variance of the preposterior mean have been estimated, $p\left(\mu_t^{\bm X}\right)$ can be approximated using \emph{moment matching} \cite{CetinkayaThiele:2016,Fengestal:2015}. In general, this involves approximating $p\left(\mu_t^{\bm X}\right)$ by taking a known distribution and ensuring that the mean and variance of this distribution equal the mean and variance of $p\left(\mu_t^{\bm X}\right)$. For example, this \emph{moment matching} method could involve approximating the distribution of the preposterior mean by a Normal distribution with the correct mean and variance.

However, this strategy is unlikely to give accurate estimates, as the EVSI is strongly influenced by the tails of the distribution of the preposterior mean because this is where the optimal decision is most likely to change. Therefore, the EVSI estimate will be significantly improved if the distribution of the preposterior mean is approximated by moment matching with the PSA samples for $\nb_t^{\bm\theta}$. As these samples represent our initial beliefs about the net benefit of each treatment and, as seen in \S\ref{Example}, our initial beliefs have a strong impact on the shape of $p(\mu_t^{\x})$. 

\subsubsection{Linear transformation to moment match}
In practice, we suggest that a linear transformation of the NB$_t^{\bm\theta}$ should be used to estimate the distribution of the preposterior mean by moment matching. This involves estimating the constants $a$ and $b$ such that $a\, \mbox{NB}_t^{\bm\theta} + b$ has the same mean and variance as $p(\mu_t^{\x})$:
\begin{align*}\mbox{E}_{\bm\theta}\left[a\ \mbox{NB}_t^{\bm\theta}+b\right]=\mbox{E}_{\bm X}\left[\mu^{\bm X}_t\right] &\Rightarrow 
a\mbox{E}_{\bm\theta}\left[\mbox{NB}_t^{\bm\theta}\right]+b = \mbox{E}_{\bm\theta}\left[\mbox{NB}_t^{\bm\theta}\right] \\
\mbox{Var}_{\bm\theta}[a\ \mbox{NB}_t^{\bm\theta}+b]=\mbox{Var}_{\bm X}[\mu^{\bm X}_t] &\Rightarrow a^2 \mbox{Var}_{\bm\theta}[\mbox{NB}_t^{\bm\theta}]= \sigma^2,\end{align*} where $\sigma^2$ is the variance of the preposterior mean that is calculated as the difference between the variance of the net benefit and the expected posterior variance. Solving for $a$ and $b$ yields
\begin{equation} 
a=\sqrt{\frac{\mbox{Var}_{\bm X}\left[\mu_t^{\bm X}\right]}{\mbox{Var}_{\bm\theta}\left[\mbox{NB}_t^{\bm\theta}\right]}} = \frac{\sigma}{\sqrt{\mbox{Var}_{\bm\theta}\left[\mbox{NB}_t^{\bm\theta}\right]}}\quad \mbox{and} \quad b=\mbox{E}_{\bm\theta}\left[\mbox{NB}_t^{\bm\theta}\right](1-a),\label{constants}
\end{equation} 
which depend on the expectation and variance of the net benefit and its  expected posterior variance.

Interestingly, these constants allow for a relatively simple interpretation of the approximation of $p\left(\mu_t^{\bm X}\right)$. The constant $a$ is related to the reduction in the net benefit variance that will be obtained by learning $\x$. This means that the more information $\x$ contains about the net benefit, the higher the value of $a$. In some ways, we can think of this as the amount of variance that is ``explained'' by the fact that we have not yet observed $\x$.

The constant $b$ is then the mean of the preposterior distribution weighted by one minus this \emph{explained variance} where the weight is related to how much information is contained in $\bm X$. Thus, the density of the preposterior mean is estimated as a convex combination of our initial beliefs about the net benefit and the mean of the net benefit (which is also the mean of the preposterior distribution). 

In general, the higher the sample size of $\x$, the more information it contains and therefore the higher the value of $a$. This implies that as the sample size in the data collection exercise increases, our approximation for $p(\mu_t^{\x})$ approaches the PSA samples for the net benefit, i.e.~an infinite future sample implies that the distribution of the preposterior mean is exactly equal to our initial beliefs about the net benefit. We explore this slightly counterintuitive result in detail in the supplementary material along with further motivation for our moment matching method. Note that, as a rule of thumb, the approximation is most accurate when the sample size of the future dataset is greater than 20, especially when the data are discrete, as in \S\ref{Example}.

As PSA simulations for the net benefit are generally available as a component of a full economic analysis \cite{BaioDawid:2011,Baio:2012,Andronisetal:2009}, these samples can be used to estimate the mean and variance of $\nb_t^{\bm\theta}$. These mean and variance values can then used to calculate $a$ and $b$. Additionally, the PSA simulations can be rescaled using these constants  so the only additional element required to approximate the distribution of the preposterior mean is an estimate of the expected variance of the posterior net benefit across different possible future samples~$\bm X$. 

\subsection{Estimating the expected variance of the posterior net benefit}\label{pre-post var}
We suggest that the expected variance of the posterior net benefit should be estimated by simulation. While this seems computationally intensive, the number of simulations can be kept fairly low as the posterior variance is relatively \emph{stable} across different future samples $\bm X$, meaning that the variance does not change substantially for different samples. This stability is most extreme when the data distribution and our initial beliefs are both normal. In this case, the posterior variance is the same irrespective of the posterior mean as it is only related to the variance of the sample $\bm X$ and the sample size. 

In non-normal settings (e.g.~involving the collection of data for costs or utility measures), the posterior variance is no longer unrelated to the mean of the sample $\x$. However, as demonstrated in the supplementary material, a small number of ``posterior'' samples (around 20-50) can be used to estimate the expected posterior variance accurately, even in highly non-normal settings. Crucially, this result is dependent on intelligently choosing which posterior samples to use.

Specifically, we need to calculate the expected posterior variance: 
\[\mbox{E}_{\bm X}\left[\mbox{Var}_{\bm\theta\mid\bm X} \left[\mbox{NB}_t^{\bm\theta}\right]\right] = \mbox{E}_{\bm\theta}\left[\mbox{E}_{\bm X \mid \bm\theta}\left[\mbox{Var}_{\bm\theta\mid\bm X} \left[\mbox{NB}_t^{\bm\theta}\right]\right]\right], \] where the two outer expectations on the right hand side estimate the expectation over all possible values of $\x$. In general, we need to choose a small number of $\bm\theta$ values whilst ensuring that the full range of $\bm\theta$ values is explored before sampling \emph{one} future sample conditional on each value of $\bm\theta$. This estimates the expected posterior variance accurately using quadrature \cite{DavisRabinowitz:2007}.

In most health economic models $\bm\theta$ contains a large number of parameters, sometimes upwards of 1000 \cite{Walenskyetal:2010}. However, most trials are only designed to gather information about a small number of these parameters (in fact, trials often focus on a single primary outcome to be estimated reliably). We denote these targeted parameters $\bm\phi$ and it is only these parameters that we must explore intelligently to estimate the expected posterior variance. Specifically, it is suggested that $Q$ evenly spaced values for $\bm\phi$ should be chosen from the PSA samples for $\bm\phi$. These $Q$ values are then used to simulate one future sample from $\bm X\mid\bm\phi$ for each $\bm\phi$ value, resulting in $Q$ simulated datasets.

A posterior distribution is then found for each simulated future sample. This gives $Q$ posterior distributions, which will usually be found using MCMC simulation through computer programs such as \texttt{BUGS} \cite{Lunnetal:2000}, \texttt{JAGS} \cite{Plummer:2003} or \texttt{Stan} \cite{Carpenteretal:2016}. For each posterior distribution, we then need to calculate the variance of the net benefit in order to approximate the EVSI using moment matching. In general, $Q$ should be greater than 20 to obtain sufficiently accurate estimates of the variance. However, as $Q$ increases, the EVSI estimate becomes more accurate, although at the cost of longer computational time. Estimating the variance in this manner implies that our method can be used irrespective of the underlying model structure, provided it is computationally feasible to find the posterior distribution by MCMC simulation $Q$ times.

\subsection{Calculating the EVSI for a specific set of treatment options}
To calculate the EVSI using this estimated distribution of the preposterior mean, it is necessary to determine the joint distribution of the preposterior mean net benefit across all the different treatment options. Theoretically, this adds little complication to the method already described in that a posterior variance-covariance matrix for all the net benefits must be estimated, instead of a scalar variance. 

However, the EVSI estimate is more stable computationally if the \textit{incremental} net benefit (INB) of treatment $r$ versus treatment $s$, defined as $\mbox{INB}^{\bm\theta}=\mbox{NB}_r^{\bm\theta}-\mbox{NB}_s^{\bm\theta}$, is used instead. In this case, the optimal treatment, under current or future information, is found by comparing the INB with 0. For example, if only two treatment options are considered, then treatment $t=2$ is optimal if the INB is positive while $t=1$ is optimal if it is negative. More importantly, if the distribution of the preposterior mean of the INB is estimated using the moment matching method, the EVSI is computed as \[\mbox{E}_{\bm X}\left[\max\left\{0,\mbox{E}_{\bm\theta\mid \bm X}\left[\mbox{INB}^{\bm\theta}\right]\right\}\right] - \max\left\{0,\mbox{E}_{\bm\theta}\left[\mbox{INB}^{\bm\theta}\right]\right\},\] where $\mu^{\bm X}= \mbox{E}_{\bm\theta\mid \bm X}\left[\mbox{INB}^{\bm\theta}\right]$ and the mean and variance of the INB are scalars rather than matrices.

When more than two treatment options are being compared, working with the INB reduces the size of the posterior variance-covariance matrix that must be estimated, leading to greater numerical stability. For example, if three treatment options are considered, using the INB implies that only three unique elements must be estimated in the posterior variance-covariance matrix rather than 6.

\subsection{Nuisance Parameters}\label{ind-para} 
Finally, an extension to the method presented so far is needed to reflect
the fact that the new research only informs a small subset of the model parameters $\bm\phi$, as discussed in \S\ref{pre-post var}. In general, it is unlikely that the distribution for INB$^{\bm\theta}$, conditional on all the underlying model parameters $\bm\theta$, will be the same as the distribution of the preposterior mean for data informing $\bm\phi$. Therefore, it is necessary to determine the distribution of the INB conditional solely on $\bm\phi$, indicated as INB$^{\bm\phi}$. It is then these samples of INB$^{\bm\phi}$ that  are rescaled linearly to moment match with the distribution of the preposterior mean. 

To demonstrate the difference between INB$^{\bm\theta}$ and INB$^{\bm\phi}$, we consider a simple two parameter model where: $\phi \sim \mbox{Beta}(1,4)$, $\psi \sim \mbox{Normal}(-0.5,1)$, $\mbox{NB}_0^{\bm\theta} = 10\,000\psi - 4\,000$; $\mbox{NB}_1^{\bm\theta} = 10\,000\phi-6\,500$; and $\mbox{INB}^{\bm\theta}~=~10\,000~(\phi~-~\psi)~-~2\,500$.
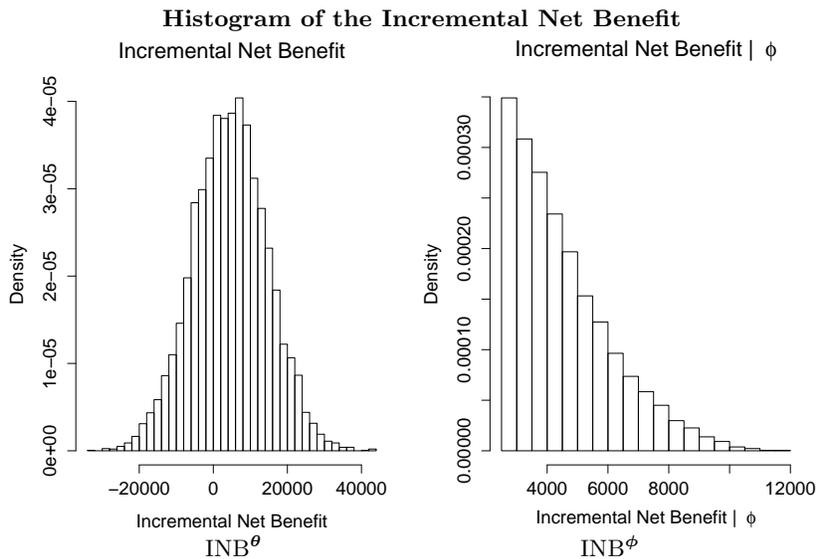
\begin{figure}[!h]
\begin{center}
\begin{tikzpicture}
\node at (0,0){\includegraphics[width=11cm]{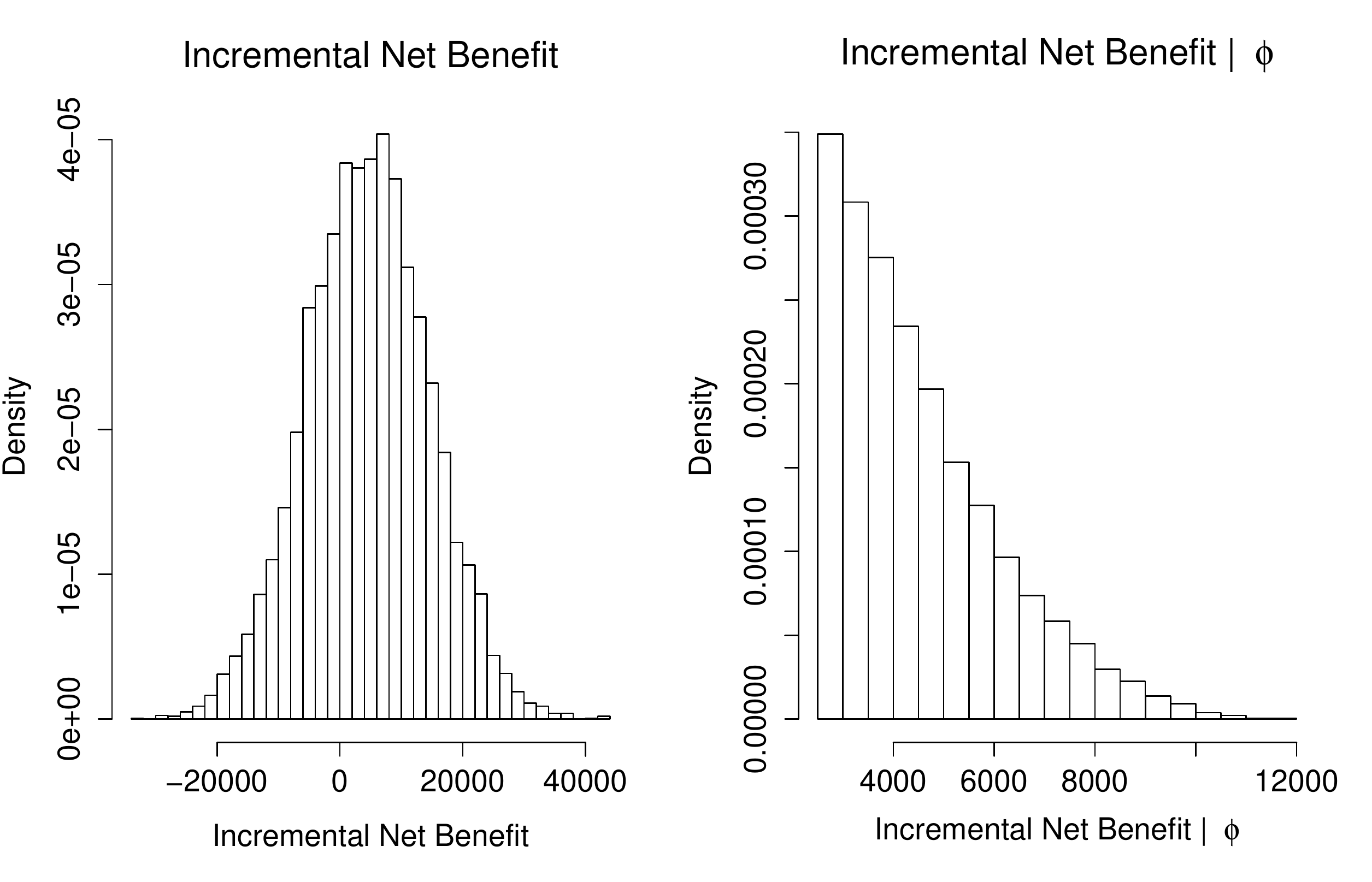}};
\node at (-2.5,-3.5){\small $\mbox{INB}^{\bm\theta}$};
\node at (2.5,-3.5){\small $\mbox{INB}^{\bm\phi}$};
\node at (0,3.5){\small \textbf{Histogram of the Incremental Net Benefit}};
\end{tikzpicture}
\end{center}
\caption{The distribution of the incremental net benefit conditional on both the model parameters $\bm\theta$ (LHS) and conditional on the parameter of interest $\phi$ (RHS).}
\label{ind-ex}
\end{figure}

Figure \ref{ind-ex} shows that the prior for INB$^{\bm\theta}$ is approximately normal --- particularly in the tails --- while, as the INB is linear in $\phi$, the distribution of INB$^{\bm\phi}$ is a shifted and scaled Beta distribution, which is non-normal. Therefore, if INB$^{\bm\theta}$ was used to approximate the distribution of the preposterior mean for new research that only informs the parameter $\phi$ then the shape of that approximate distribution would be incorrect, leading to inaccurate estimates for the EVSI.

To identify a more appropriate shape for the distribution of the preposterior mean in the presence of nuisance parameters, i.e.~the parameters $\bm\psi$ such that $\bm\theta = (\bm\phi,\bm\psi)$, the uncertainty due to $\bm\psi$ should be marginalised out by calculating:
\begin{equation}
\mbox{INB}^{\bm\phi}=\mbox{E}_{\bm\psi\mid\bm\phi}\left[\mbox{INB}^{\bm\theta}\right].\label{marginalise}
\end{equation}
While calculating this expectation can be computationally intensive, it should be estimated before proceeding to calculate the EVSI \cite{Tuffahaetal:2016}. This is because calculating the EVSI requires a study design for the future research which is time consuming to determine for a realistic trial. Therefore, before designing the trial, the value of resolving \emph{all} the uncertainty in $\bm\phi$ should be calculated. This is known as the Expected Value of \emph{Partial} Perfect Information (EVPPI) and is based on the expectation in \eqref{marginalise}. Consequently, samples of INB$^{\bm\phi}$ should already have been calculated. However, if these values are not available, Strong \textit{et al.} \cite{StrongOakley:2014} or Heath \textit{et al.} \cite{Heathetal:2016} offer computationally efficient procedures for estimating this expectation using non-parametric regression, with code available in the \texttt{R} package \texttt{bcea} \cite{BCEA:2017} or as a stand-alone function \cite{Strong:2012:Code}. 

Finally, to rescale the sample of the INB$^{\bm\phi}$, rather than INB$^{\bm\theta}$, the constant $a$ from \S\ref{constants} becomes: \[a = \frac{\sigma}{\sqrt{\mbox{Var}_{\bm\phi}\mbox{INB}^{\bm\phi}}},\] where $\sigma$ is the variance of the preposterior mean. Note that $\sigma$ must still be calculated as a function of the variance of INB$^{\bm\theta}$ --- not the variance of INB$^{\bm\phi}$ --- and the posterior variance, as outlined in \S\ref{pre-post var}.

\section{Case Study: Ades et al. Decision Tree Model}\label{toy}
To demonstrate the effectiveness of our methodology, we use a decision tree model developed in Ades \textit{et al.} \cite{Adesetal:2004}. This model has two treatment options, a standard of care and a new treatment, aimed at avoiding a critical event. This critical event leads to a reduction in QALYs for the remainder of the patient's life. The new treatment reduces the probability of the critical event but the patient may also experience side effects which give a short term reduction in QALYs along with a direct cost of additional treatment. The model has 11 parameters, of which 4 are subject to uncertainty which is then modelled using 4 mutually independent distributions; a complete model description is given in Ades \textit{et al.} \cite{Adesetal:2004} or Strong \textit{et al.} \cite{Strongetal:2015}.

For this case study, we consider four different data collection exercises, the first three have been tackled by Ades \textit{et al.} and the fourth investigates the moment matching method when $\bm\phi$ is not unidimensional;
\begin{enumerate}
\item To reduce uncertainty in the probability of side effects for the new treatment $\phi_1$, 60 patients are given the new treatment and the number who suffer from side effects is recorded.
\item To reduce uncertainty in the quality of life after the critical event $\phi_2$, the quality of life for 100 patients who experienced the event is recorded.
\item To reduce uncertainty in the odds ratio of the effectiveness of the two treatments $\phi_3$, a randomised control trial with 200 patients on each arm is undertaken.
\item To reduce uncertainty in the probability of the critical events in both treatment arms $\phi_3^C,\phi_3^T$, the same randomised control trial is undertaken but the analysis informs these two probabilities directly.
\end{enumerate}
For a full description of the distributional assumptions for these studies, particularly the difference between study 3 and study 4, see Strong \textit{et al.}~ \cite{Strongetal:2015}.

\subsection{Computations}

To calculate the EVSI using the moment matching methodology, 1\,000\,000 simulations were taken for the 4 stochastic model parameters. These were then combined with the other seven parameters to calculate 1\,000\,000 simulations for the INB under current information. These simulations for INB$^{\bm\theta}$ were used to find INB$^{\bm\phi}$ using GAM regression \cite{HastieTibshirani:1990} obtained with the \texttt{gam} function from the \texttt{mgcv} package \cite{mgcv:2016} in \texttt{R}. The simulations were also used to find the mean and variance of INB$^{\bm\theta}$.

To estimate the preposterior variance, the expected posterior variance was estimated using MCMC procedures with $Q$, the number of $\bm\phi$ values, equal to 30. This was achieved using \texttt{JAGS} through \texttt{R} \cite{Plummer:rjags} with 10\,000 simulations from each posterior distribution and 1\,000 simulations used as burn-in. This means that in total 1\,330\,000 simulations were used to estimate the EVSI in this example, although note that the PSA simulations were reused for each EVSI calculation. Fewer simulations could be used but as this example has a small computational cost it was possible to use this number of simulations to improve the accuracy of the methods. The computational time required to estimate the posterior variance with this number of simulations was between 3.8-6.1 seconds. 

To assess the accuracy of our method, the expected posterior variance was estimated using the above procedure 1\,000 times for each trial. Each of these estimates for the expected posterior variance was then used to approximate the EVSI using moment matching to give a distribution for the EVSI estimate obtained using moment matching. In a standard analysis, therefore, our method would produce a point estimate for the EVSI rather than the distributions given in \S\ref{results}.

To determine the accuracy of our method, we compared with the computationally intensive two-step nested Monte Carlo procedure \cite{Brennanetal:2007} and the Strong \textit{et al.} method \cite{Strongetal:2015} based on sufficient statistics and non-parametric regression which is an accurate and efficient estimation method for models where the data can be summarised using a low dimensional sufficient statistic. For the experiments targeting $\phi_1$, $\phi_2$ and $\phi_3$, the comparator values are taken directly from the Strong \textit{et al.}~paper where the estimates are based on $10^{10}$ and $10^6$ simulations respectively. For the two parameter EVSI estimate, the results for these two methods were obtained using the same number of simulations using code given in the supplementary material. The computational times to obtain these estimates were 8.4 seconds for the Strong \textit{et al.}~method and 207\,595 seconds (approximately 2.4 days) for the nested Monte Carlo method.

\subsection{Results}\label{results}
\begin{figure}[!h]
\includegraphics[width=15cm]{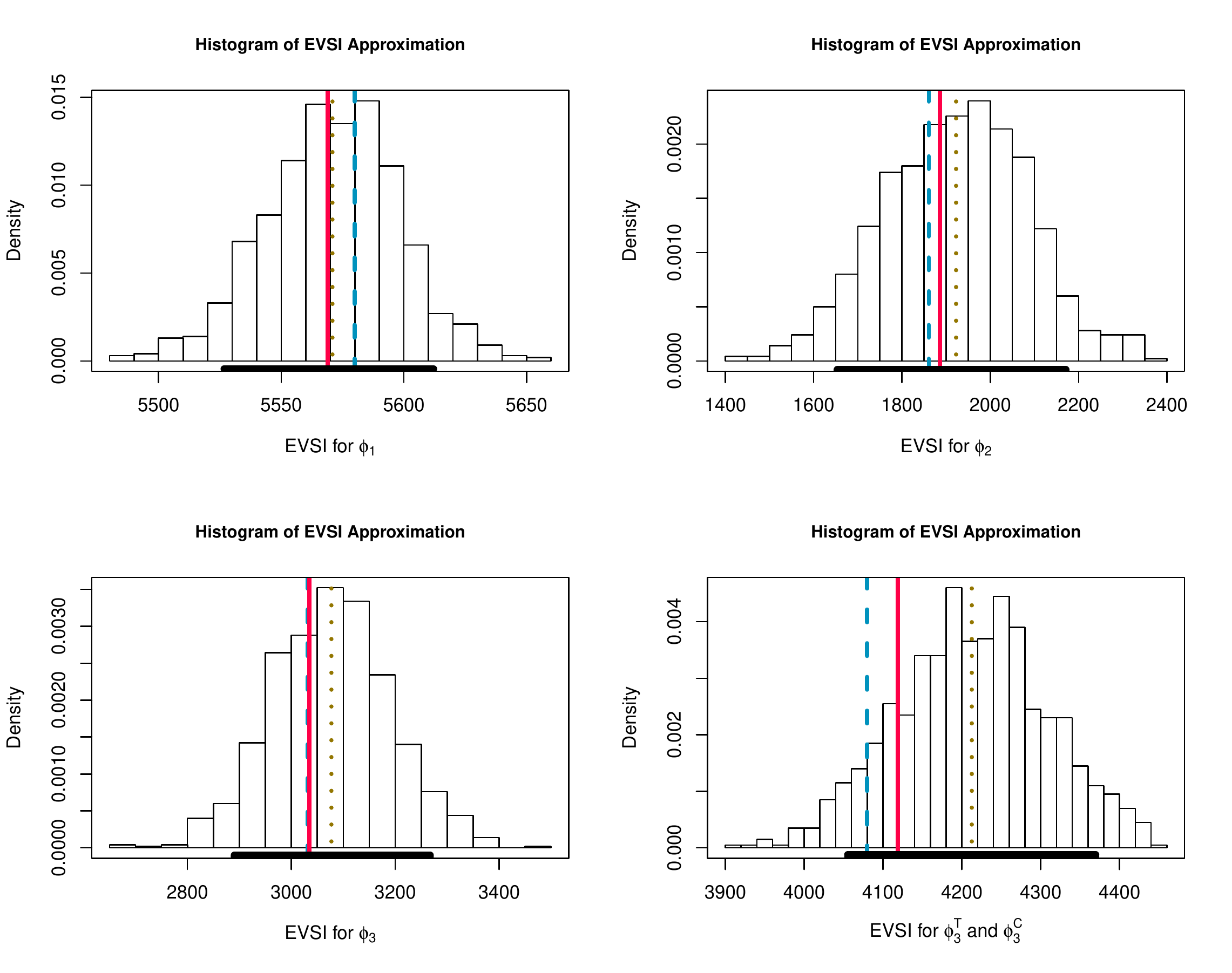}
\caption{The ``sampling distribution'' of the EVSI conditional on the distribution over the different estimates for the variance of the preposterior mean for $\phi_1$, $\phi_2$ and $\phi_3$ for the Ades \textit{et al.} example \cite{Adesetal:2004}. The solid red line represents the EVSI calculated using Monte Carlo methods and $10^{10}$ simulations. The dashed blue line represents the EVSI estimate obtained using the Strong \textit{et al.}~method with $10^6$ simulations. The dotted brown line represents the average EVSI estimate obtained using moment matching and the solid black horizontal line represents the 90\% interval for this estimation method. The comparator methods (MC and Strong \textit{et al.}) are taken from Strong \textit{et al.}~\cite{Strongetal:2015} for the first three graphics.}
\label{est-Ades}
\end{figure}

Figure~\ref{est-Ades} plots the sampling distribution of the EVSI estimate, over the different estimates of the expected posterior variance, obtained using our moment matching method for the four alternative study designs. The solid line gives the value of the estimate obtained by nested Monte Carlo (which can be considered as the ``truth''), the dashed line is the estimate obtained using the Strong \textit{et al.}~method and the dotted line is the average moment matching estimate. Evidently, the moment matching method is in line with these two alternative estimation methods for the studies considered, although the two-parameter estimation has a slight upward bias. Nonetheless, the true EVSI value is well within the 90\% interval, represented by the solid black line on the axis in Figure \ref{est-Ades}.

Figure~\ref{est-Ades} also demonstrates that the EVSI estimate calculated using moment matching is less accurate for experiments with smaller EVSI --- as the EVSI gets smaller, the 90\% intervals get wider. This is because the estimate is based on the difference between the prior variance and the expected posterior variance. When the EVSI is small then this difference is also small and therefore the posterior and prior variance needs to be estimated with greater precision because the difference can be greatly affected by the Monte Carlo error in the variance estimation. It is important to note that the accuracy demonstrated in Figure~\ref{est-Ades} depends on accurately estimating the prior variance. Therefore, our moment matching method should not be used if the initial PSA simulation size is very small.

\section{Discussion}
We have presented an estimation method that reduces the significant computational burden required to estimate the EVSI accurately. This method is based on using information already available to the researcher, typically in the form of PSA simulations, a relationship that is also exploited by Menzies \cite{Menzies:2016} despite large differences between the two methods. Our method involves moment matching by performing a linear transformation of the PSA simulations for the incremental net benefit. To perform this matching, nested sampling is used to calculate the expected posterior variance across different samples. For the example in this paper, the required number of nested simulations was reduced from $10^6$ to 30, offering a significant computational saving especially in more complex models where running the model has a computational cost. 

\subsection{When should this method be used?}
Several methods have been proposed for estimating the EVSI. A large number of these are restricted to models that fulfil certain conditions, either on the economic model being used or the structure of the data, or potentially both \cite{Weltonetal:2014,BrennanKharroubi:2007b,Adesetal:2004,Jalaletal:2015}. The moment matching method is likely to have a higher computational cost than these methods but can be used in any setting. Current general purpose methods, i.e.~those that can be used irrespective of the underlying model, include the one presented by Strong \textit{et al.}, used as a comparator in this paper, the one developed by Menzies \cite{Menzies:2016} and a method based on Laplace approximations to bypass the use of MCMC to perform the Bayesian updating \cite{Kharroubietal:2011,BrennanKharroubi:2007}.

Firstly, provided a low-dimensional sufficient statistic is available to capture all the information contained in the sample, the Strong \textit{et al.}~\cite{Strongetal:2015} method typically outperforms the other methods in terms of computational time. However, it may be challenging to determine an appropriate sufficient statistic, particularly for complex designs. The Menzies \cite{Menzies:2016} method also avoids the need for additional model runs but does involve a large number of matrix operations, which has a significant computational cost when the number of PSA simulations is large. Therefore, the moment matching method can be significantly faster than this alternative when the model itself has a relatively small computational cost. Conversely, if the underlying health economic model is computationally intensive to run it would not feasible to rerun the model a large number of times and the Menzies method could be used.

Finally, using Laplace approximations requires $T(2p+1)$ model evaluations, where $p$ is the number of model parameters and $T$ the number of treatment options \cite{BrennanKharroubi:2007}. It also involves numerical optimisation, which has an associated computational cost. Nonetheless, this method could be more efficient than our moment matching method in health economic models with a small number of underlying parameters but a long run~time. 

In conclusion, the moment matching method presented in this paper is a computationally efficient method for estimating the EVSI and is most helpful in cases where a sufficient statistic for the data is difficult to specify but the computational time taken to run the model itself is relatively small, although it can be competitive in terms of computational time for simple models. Therefore, we believe this method is an important addition to EVSI literature as it is based solely on Bayesian updating which typically must be designed in order to analyse the data once the experiment has been performed. If the research goes ahead, the calculation of the EVSI will aid in modelling trial results as researchers can use the procedure that has already been defined.

\section{Acknowledgements}
Financial support for this study was provided in part by grants from the UK Engineering and Physical Sciences Research Council (EPSRC) [Anna Heath] and Mapi [Dr. Gianluca Baio] thorough an unrestricted research grant at University College London. The funding agreement ensured the authors' independence in designing the study, interpreting the data, writing, and publishing the report. The authors would also like to thank the reviewers for their invaluable suggestions.

\bibliographystyle{plain}
\bibliography{bib}

\appendix
\section{The distribution of the preposterior mean}
To explore the similarity between the prior and the distribution of the preposterior mean two examples are presented.
\subsection{Exponential Gamma Model}

\begin{figure}[!h]
\includegraphics[width=13cm]{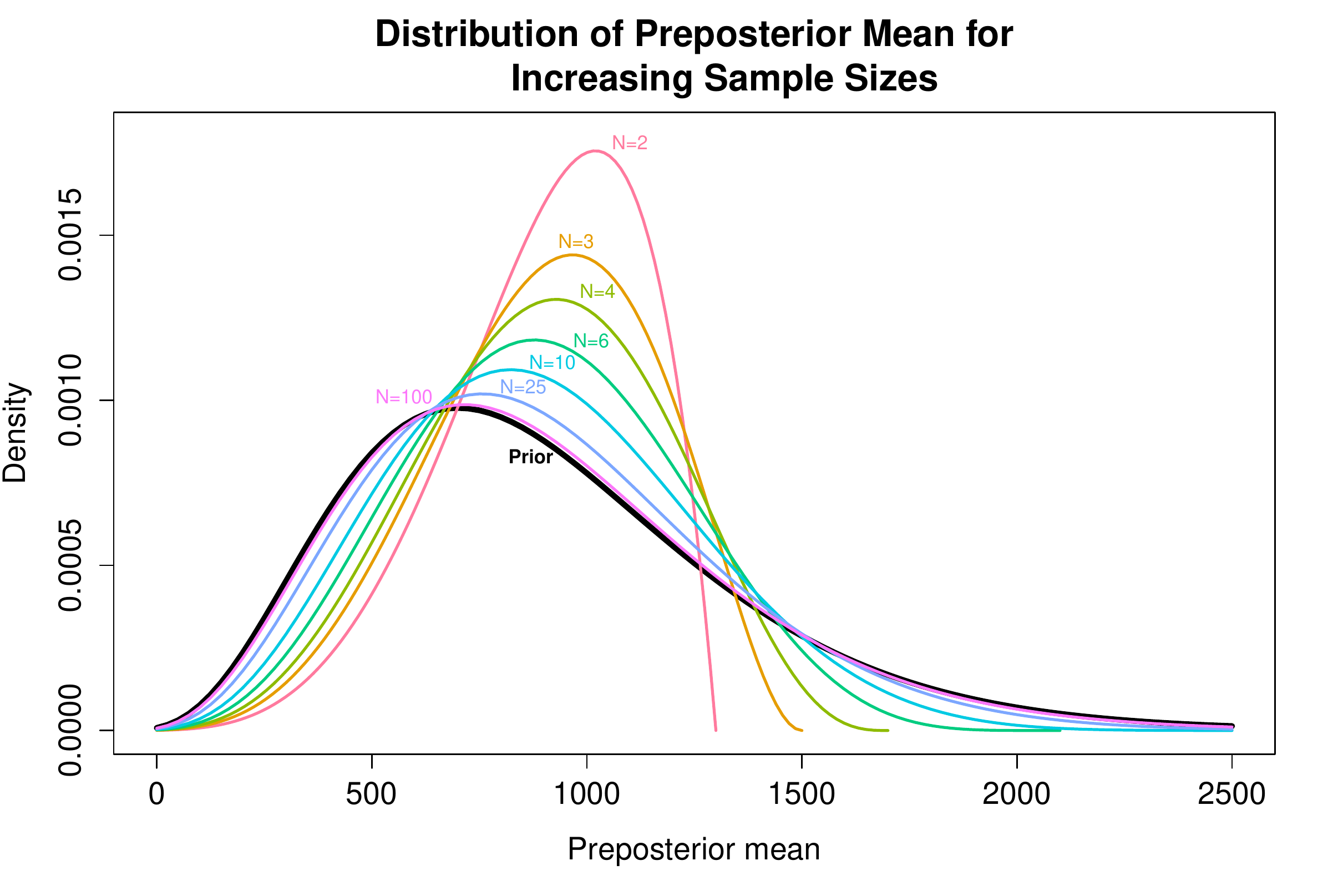}
\caption{The distribution of the exact preposterior mean for different samples sizes using Exponential-Gamma conjugacy, with the prior for the net benefit marked in black.}
\label{prepostmean}
\end{figure}

For the first example, a Gamma prior is assumed for the parameter of interest $\theta~\sim~Gamma(\alpha,\beta)$. The data collection exercise is then assumed to be $N$ independent observations from an exponential distribution conditional on $\theta$; $X_j~\sim~Exp(\theta)$ with $j=1,\dots,N$. The distribution of the preposterior mean is then considered for different values of $N$, where the two net benefit functions are: \[\mbox{NB}^\theta_0=c_0 \quad \mbox{and} \quad \mbox{NB}^\theta_1=k\theta-c_1.\] Figure \ref{prepostmean} presents the distribution of the preposterior mean for NB$_1^\theta$ for $\alpha=5$, $\beta=1$, $k=200$, $c_0=900$ and $c_1=100$.

Firstly, observe that the distribution of the preposterior mean approaches the distribution of the prior for NB$_1$ the sample size of the future data increases. Additionally, notice that the distribution become \emph{more} variable as the sample size increases. These two comments are seemingly at odds with our intuition about the distribution of a mean --- that the distribution becomes more normal and less variable as the sample size of the data increases. 

However, these counterintuitive properties hold because the ``strength'' of the data increases as the sample size increases. This implies that the posterior mean can deviate further from the prior mean as the sample size increases. The distribution of the preposterior mean is the distribution over all possible future means, which can be further from the prior mean as the size of the future data collection increases. More importantly, the distribution of the preposterior mean tends to the prior because an infinite data collection exercise would determine the exact value of the parameter $\theta$. However, at this current state of knowledge, i.e.~\emph{before} the data collection has taken place, we believe that $\theta$ follows the prior distribution. Therefore, the distribution of the mean with this infinite data collection exercise is equal to the prior. This can be seen analytically using normal-normal conjugacy.

\subsection{Normal Normal Model}
For this example, the prior for $\theta$ is taken as \[\theta \sim N(\theta_0,\sigma^2_\theta)\] while the data collection exercise is $N$ independent samples from \[X_i \sim N(\theta,\sigma_{\bm X}^2).\] It is then assumed that the variances  $\sigma^2_\theta$ and $\sigma^2_{\bm X}$ are known. This implies that the sample mean of $\bm X$ has a normal distribution, conditional $\theta$ and that the prior-predictive distribution for $\bar{\bm X}$ is \[\bar{\bm X} \sim N\left(\theta_0,\sigma^2_\theta+\frac{\sigma^2_{\bm X}}{N}\right).\]
Finally, the net benefit functions for this example are given by \[\mbox{NB}_0^{\theta} = 0 \quad \mbox{and} \quad \mbox{NB}_1^{\theta} = k\theta-c.\]

The preposterior mean for the INB is then\[\mbox{E}_{\theta\mid\bm X}(\mbox{INB}^\theta) = k\left(\frac{\sigma^2_{\bm X}}{\sigma^2_{\bm X}+N\sigma^2_{\theta}} \theta_0 +\frac{\sigma^2_{\theta}}{\frac{\sigma^2_{\bm X}}{N}+\sigma^2_{\theta}}\bar{\bm X}\right) - c,\] which is a linear function of a normal distribution. Therefore, the distribution of the preposterior mean is normal with mean and variance equal to the mean and variance of the preposterior mean;
\[\mbox{E}_{\theta\mid\bm X}(\mbox{INB}^\theta) \sim N\left( k\theta_0-c, k^2\frac{\sigma_\theta^4}{\frac{\sigma_{\bm X}^2}{N}+\sigma_\theta^2}\right).\]

Firstly, note that the moment matching approximation for the distribution of the preposterior mean is $a\,\mbox{INB}^{\theta}+b$. This is also a linear combination of a normal distribution implying that this approximation is exactly equal to the true distribution of the preposterior mean for the normal-normal conjugate setting. 

Secondly, note that as $N\to\infty$, it is clear that the variance of the preposterior mean INB$^\theta$ tends to $k^2\sigma_\theta^2$, meaning that the distribution of the preposterior mean tends to the prior for the incremental net benefit. Another way to think of this is that, as $N\to\infty$, the sample mean $\bar{\bm X}$ is tends to the mean of $X_i$. However, in a Bayesian setting, the underlying mean $\theta$ is subject to uncertainty modeled using the prior. 

Additionally, observe that, as the sample size increases, the denominator of the variance decreases. Clearly, therefore, the distribution of the preposterior mean gets \emph{more} variable as the sample size increases. Again, this confirms as more information is contained in the data, i.e.~the sample size increases, the posterior mean can be ``pulled'' further from the prior mean and so the distribution of the preposterior mean becomes more variable.

Finally, it is trivial to see how the distribution of the preposterior mean is dependent on the prior. Firstly, is it centred on the prior mean. In addition to this, the variance of the distribution of the preposterior mean is strongly influenced by the prior variance for realistic sample sizes. Finally, the distributional assumptions for the prior clearly impact the distribution of the preposterior mean. Therefore, the moment matching method simply utilises this information as typically we will already have access to samples from the prior. 

These analytical results give confirmation that the moment matching method will be suitable in normal-normal settings. This next section briefly discusses why this method is suitable in other settings. 

\subsection{Using moment matching to estimate the distribution of the preposterior mean}\label{proof}
Firstly, consider whether the approximation is suitable in two extreme cases. At one extreme, assume that $\bm X$ is independent of the underlying model parameters: $p(\bm X\mid\bm\theta)=p(\bm X)$. Evidently, this setting would never occur as decision makers only consider data collection that would aid the decision making process. Nevertheless, if the sample is independent of the model parameters then the distribution of the preposterior mean is a point mass at the prior mean \[\mbox{E}_{\bm\theta\mid\bm X}\left[\mbox{NB}_t^{\bm\theta}\right]=\mbox{E}_{\bm\theta}\left[\mbox{NB}_t^{\bm\theta}\right],\] by the condition of independence.

Using the definition for the constants $a$ and $b$ from the main paper, note that \[a=\sqrt{\frac{\mbox{Var}\left[\mbox{E}_{\bm\theta\mid\bm X}\left[\mbox{NB}_t^{\bm\theta}\right]\right]}{\mbox{Var}_{\bm\theta}\left[\mbox{NB}_t^{\bm\theta}\right]}}=0\]
and
\[b=\mbox{E}_{\bm\theta}\left[\mbox{NB}_t^{\bm\theta}\right](1-a)=\mbox{E}_{\bm\theta}\left[\mbox{NB}_t^{\bm\theta}\right].\]
Therefore, the moment matching approximation is also equal to the prior mean. This means that, if $\bm X$ is independent, the moment matching approximation is exact. Clearly, this is a relatively unimportant result in practice but it does indicate that the approximation is roughly accurate when the variance of the preposterior mean is small.

At the other end of the scale, it can be shown that  the moment matching approximation is exact when the sample is deterministically linked to the model parameters, i.e.~$\bm X=h(\bm\theta)$ for some $h(\cdot)$. In this setting the \emph{conditional} mean for the net benefit is equal to the net benefit since, if the value for $\bm X$ is known, then the exact NB$_t^{\bm\theta}$ value is also known. In a similar manner to above it can be shown that $a=1$ and $b=0$ so the approximation for the distribution of the preposterior mean is equal to NB$_t^{\bm\theta}$ which, again, is exactly as required. Therefore, as the variance of the preposterior mean increases, the approximation becomes exact for all distributions. This has a practical implication since, provided the posterior is consistent, this approximation is accurate for large sample sizes.

To extend these ideas, this approximation is accurate for moderate $N$ when the posterior mean net benefit is a weighted average between the prior mean and a data summary \[\mbox{E}_{\bm\theta\mid\bm X}\left[\mbox{NB}_t^{\bm\theta}\right]=c\ \mbox{E}_{\bm\theta}\left[\mbox{NB}_t^{\bm\theta}\right]+ d\ g(\bm X),\] where $c$ and $d$ are constants and $g(\cdot)$ is an arbitrarily complex function of the data which must have a similar density to NB$_t^{\bm\theta}$. As seen previously, this is true in approximately normal settings and additionally, the first condition is true for all conjugate settings in the exponential family \cite{DiaconisYlvisaker:1979}. Determining whether the second condition is true is more challenging but in general it holds sufficiently well provided $N>20$.

The following examples demonstrate that for very small sample sizes the moment matching approximation can give biased estimates in non-normal settings. However, this bias is minimal for realistic sample sizes and decreases as the sample size $N$ increases past about 20.

\section{Examples}\label{suitable}
\subsection{Discrete Sampling}
To consider the moment matching approximating for discrete samples, the example from the paper is considered. Recall that the parameter $\theta$ is modelled using a uniform prior $\theta \sim Uniform(1,1)$ and the data have a binomial distribution $X\mid \theta \sim Bin(N,\theta)$. The two net benefit functions are then $\mbox{NB}_0^\theta = 0$ and $\mbox{NB}_1^\theta = k\theta-c$. 

In this setting, the approximation of the distribution of the preposterior mean could be poor, as $X$ is discrete. This implies that the distribution of the preposterior mean is discrete, while the prior for NB$_1^\theta$ is continuous. For example, when the binomial sample size $N=1$, there are 2 equally likely possible samples, $X=0$ and $X=1$, implying that there are two equally likely possible preposterior means; $\mu_1^0=\frac{k}{3}-c$ or $\mu_1^1=\frac{2k}{3}-c$. Clearly, this distribution can never be well approximated by a shifted and rescaled beta distribution.

Nonetheless, to investigate when a continuous approximation is suitably accurate, the EVSI is estimated for different binomial sample sizes $N$. As conjugate models are used, it is possible to calculate both the EVSI and the variance of the preposterior mean analytically. The true variance of the preposterior mean is then used to rescale the simulated observations from INB$^{\theta}$. While the variance of the preposterior mean is known, the moment matching estimator for the EVSI is still random as simulated values are taken from the prior for $\theta$. Therefore, 10\,000 different simulations of size 10\,000 were taken from the prior for $\theta$ and used to calculate the EVSI. This gives the sampling distribution of the EVSI estimator which should be centred on the true value for the EVSI.
\begin{figure}[!h]
\includegraphics[width=9.5cm]{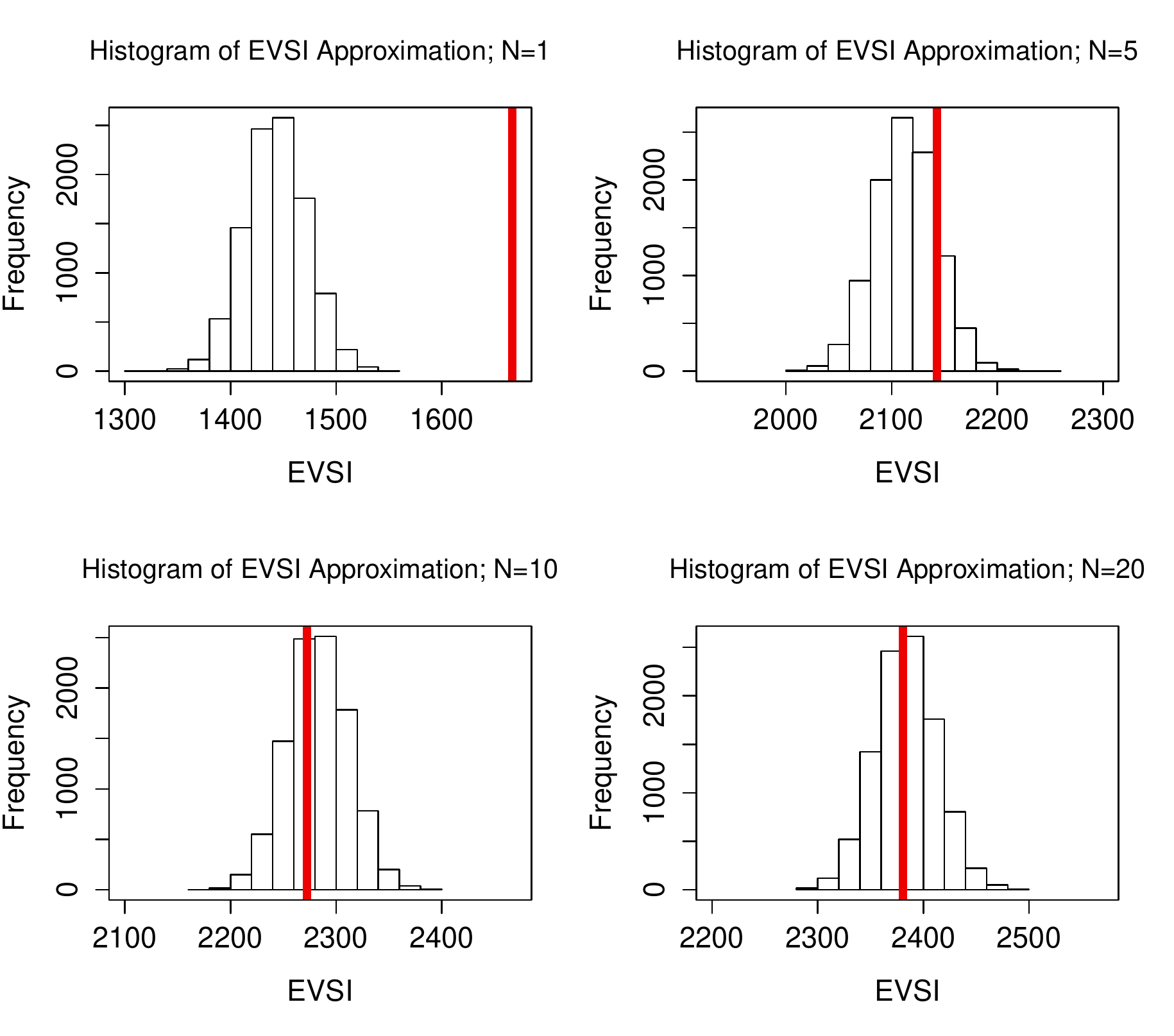}
\caption{The distribution of the EVSI estimator over 10\,000 different simulations from the prior of $\theta$ for 4 different sample sizes for $\bm X$ for the Beta-Binomial conjugate model. The red line represents the analytical value of the EVSI.}
\label{beta-bin-1}
\end{figure}

Figure \ref{beta-bin-1} shows the sampling distribution of the EVSI estimator for different $N$ where the red line gives the true EVSI. Clearly the EVSI estimator for $N=1$ has a significant downward bias as the sampling distribution does not include the true EVSI value of $1667$ (top LHS). Therefore, the weighted prior distribution is not a suitable approximation for the distribution of the preposterior mean for $N=1$. However, as $N$ increases, the bias decreases, becoming negligible for $N=10$. Therefore, even if the distribution of the preposterior mean is discrete, the moment matching approximation is suitably accurate when $N$ is sufficiently large. 

\subsection{Non-linear mean function with Exponential-Gamma Example}\label{Exp-Gamma}
To investigate the moment matching further consider, the Gamma-Exponential example: $\theta \sim Gamma(\alpha,\beta)$, $ X_j \sim Exp(\theta)$, $j=1,\dots,N$ and net benefit functions are: \[\mbox{NB}^\theta_0=c_0 \quad \mbox{and} \quad \mbox{NB}^\theta_1=k\theta-c_1,\] with $\alpha=5$, $\beta=1$, $k=200$ and $c_0=900$ and $c_1=100$ as in Figure \ref{prepostmean}. This implies that the preposterior mean is equal to \[\mu_1^{\bm X} = \mbox{E}_{\theta\mid\bm X}\left[\mbox{NB}_1^\theta\right] = k\frac{\alpha+N}{\beta+\sum_{i=1}^{N} X_i} -c_1,\] which, in turn, means that both the variance of the preposterior mean and the EVSI can be found analytically. Therefore, poor estimation of the EVSI using moment matching is because the distribution of the data summary is not similar enough to the prior.
\begin{figure}[!h]
\includegraphics[width=9.5cm]{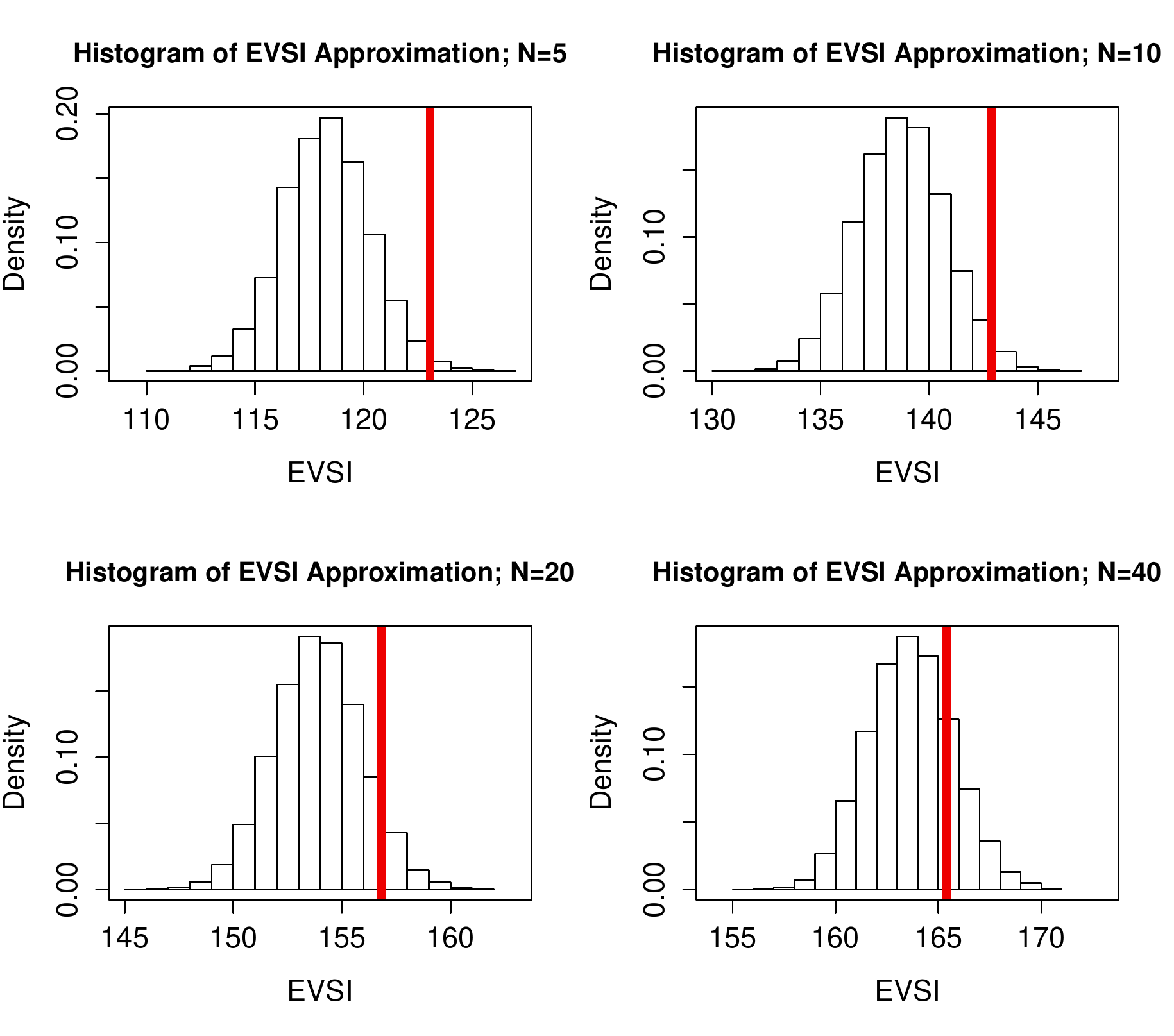}
\caption{The distribution of the EVSI estimator over 10\,000 different simulations from the prior of $\theta$ for 4 different samples sizes for $\bm X$ using the Exponential-Gamma conjugate model. The red line represents the analytic value of the EVSI,}
\label{Poisson-Gamma}
\end{figure}

Figure \ref{Poisson-Gamma} shows the sampling distribution of the EVSI values, over different prior samples, for different values of $N$. Clearly, the moment matching approximation gives a biased EVSI estimate for small samples. However, this bias is at most 4\% of the total EVSI value, meaning that it is likely to be sufficiently accurate for decision making. However, note that the EVSI estimate is slightly biased for small sample sizes so this method performs better with more realistic sample sizes. This is because the distribution of the preposterior mean tends to the prior as the sample size increases, see Figure \ref{prepostmean}.

\section{Estimating the variance of the preposterior mean}\label{non norm var}
To investigate whether quadrature is suitable when the prior for the INB$^{\bm\theta}$ is highly non-normal, we introduce a new model. In this setting, INB$^\theta= \theta^2 -5$ where $\theta$ is normal \textit{a priori} with mean 0 and precision $0.2$: $\theta \sim N\left(0,5\right)$. The data collection is then assumed to be $10$ independent observations $X_j \sim N(\theta, 1)$ for $j=1,\dots,10$. We can easily find the posterior for $\theta$ efficiently but the form of the INB induces a highly non-normal prior. In this setting, the EVSI and the variance of the preposterior mean are estimated using efficient Monte Carlo methods \cite{Adesetal:2004} with 10\,000 samples from the prior for the INB$^\theta$. The EVSI is estimated as 2.00 and the variance of the preposterior mean is 35.20. 

Recall that the estimation method for the variance of the preposterior mean requires $Q$ points spaced throughout the PSA values of $\theta$. In this practise, these are the $Q$ quantiles for $\theta$, i.e.~the $S\frac{q}{Q+1}$--th $\theta$ values in an ordered sample, with $q=1,\dots, Q$. 
\begin{figure}[!h]
\includegraphics[width=8.5cm]{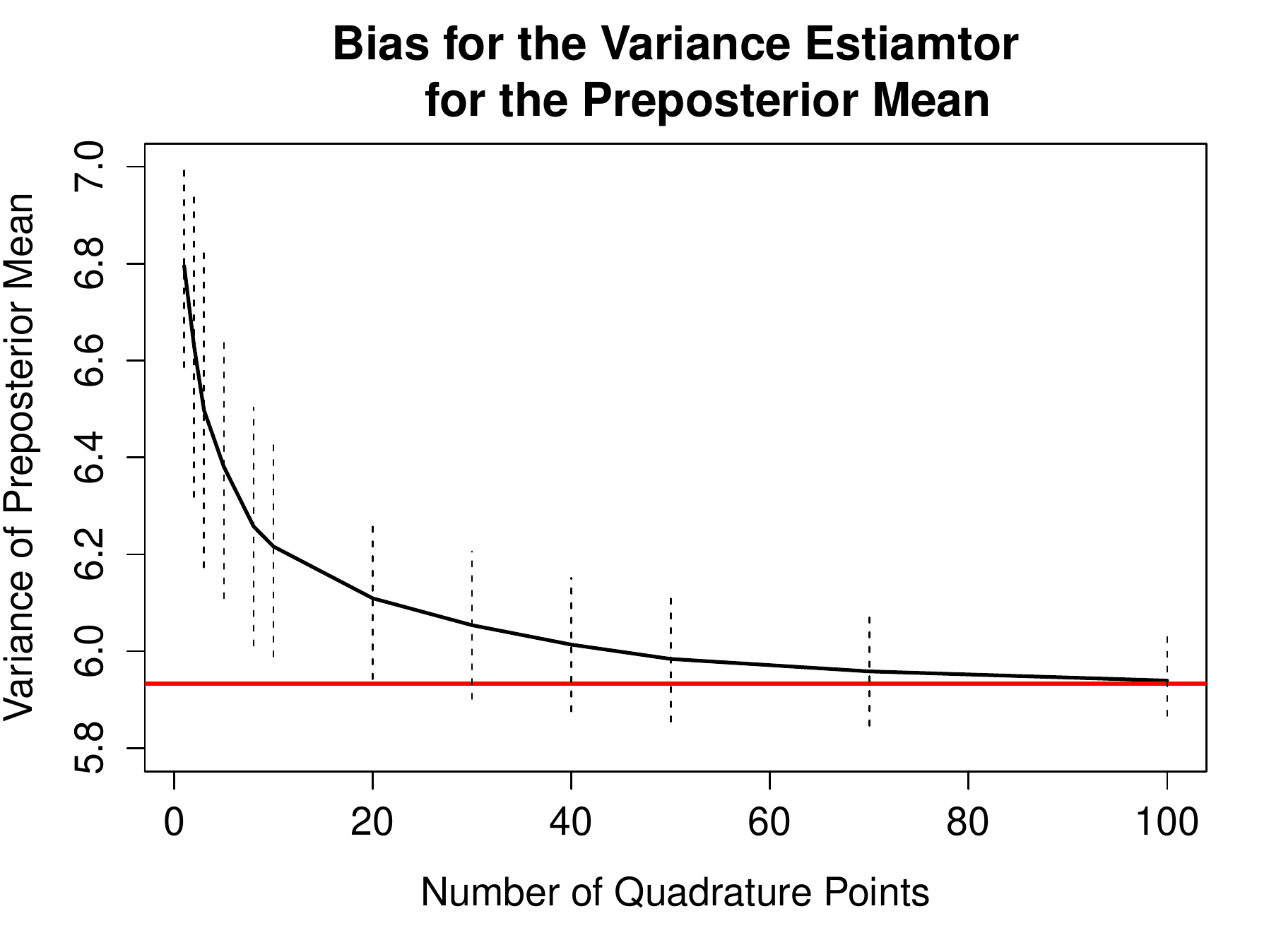}
\caption{The estimate of the variance of the preposterior mean for increasing numbers of quadrature points. The red line gives the variance of the preposterior mean calculated using all the samples in the prior for $\theta$. The dashed lines are the standard errors for the estimates of the variance of the preposterior mean.}
\label{var-samp}
\end{figure}
Figure \ref{var-samp} then shows the average estimate, over 500 posterior simulations of 1\,000, of the preposterior variance for increasing values of $Q$ up to $Q=100$ --- this means that 1\,000 simulations were taken from $Q$ different posterior distributions to estimate the variance of the preposterior mean and then this was repeated 500 times to give a measure of accuracy. The red line in Figure \ref{var-samp} is the variance of the preposterior distribution estimated using the method from Ades \textit{et al.} \cite{Adesetal:2004}. The dashed lines indicate plus or minus one standard deviation from the mean estimate of the variance of the preposterior distribution for the different values of $Q$.

In general, the estimation method for the variance of the preposterior mean produces biased estimates for small values of $Q$. However, when $Q>30$, the true variance is within one standard deviation of the average estimate for the variance of the preposterior mean. Note also that the estimate becomes more accurate as $Q$ increases but this is at the cost of computational time. Therefore, it is recommended to keep $Q$ between 30 and 50.
\begin{table}[!h]
\caption{The EVSI estimate for different numbers of posterior samples using the moment matching method.}
\begin{tabular}{|l|c|c|c|c|c|c||c|}

\hline
\centering Number of simulations& 1&2&3&5&8&10&$\infty$ \\
\hline
\centering Estimate of EVSI &2.30 &2.24 &2.20& 2.16& 2.12 &2.10&2.00\\
\hline
\centering Percentage Bias &0.15 & 0.12& 0.10 &0.08& 0.06& 0.05 &0.00\\
\hline
\hline
\centering Number of simulations&20&30&40&50 &75&100&$\infty$ \\
\hline
\centering Estimate of EVSI & 2.06& 2.05& 2.03& 2.02& 2.01& 2.01 &2.00\\
\hline
\centering Percentage Bias& 0.03 &0.02 &0.02 &0.01& 0.01& 0.01&0.00\\ \hline
\end{tabular}

\label{estimate}

\end{table}

To confirm this, table \ref{estimate} records the EVSI estimate and its bias due to using the moment matching method and the variance estimates in Figure \ref{var-samp}. Note that, the EVSI is upwardly biased by this bias drops below $0.02\%$  for $Q>30$. Therefore, a relatively small number of quadrature points can be used to estimate the variance of the preposterior mean, even in significantly non-normal settings.

\section{Code to calculate the EVSI for $\phi_3^T$ and $\phi_3^C$}

\begin{verbatim}
set.seed(2000)
library(mgcv)
library(boot)
library(R2jags)
###Prior for INB###
N<-1000000
L<-30
Qe<-inv.logit(rnorm(N,0.6,sqrt(1/6)))
Qse<-1
Ce<-200000
Ct<-15000
Cse<-100000
Pc<-rbeta(N,15,85)
Pse<-rbeta(N,3,9)
OR<-exp(rnorm(N,-1.5,sqrt(1/3)))
Pt<-inv.logit(log(Pc/(1-Pc))+log(OR))
lambda<-75000

NB1<-Pc*(lambda*L*(1+Qe)/2-Ce)+
  (1-Pc)*lambda*L

NB2<-Pse*Pt*(lambda*(L*(1+Qe)/2-Qse)-(Ct+Cse+Ce))+
  Pse*(1-Pt)*(lambda*(L-Qse)-(Ct+Cse))+
  (1-Pse)*Pt*(lambda*L*(1+Qe)/2-(Ct+Ce))+
  (1-Pse)*(1-Pt)*(lambda*L-Ct)
INB<-NB2-NB1
###Prior Variance###
pr.var<-var(INB)

###Strong et al. Method
Dc<-rbinom(1e+06,200,Pc)
Dt<-rbinom(1e+06,200,Pt)
fitted.X<-gam(INB~te(Dc,Dt))$fitted

EVSI.Strong<-mean(pmax(fitted.X,0))-max(mean(fitted.X))

###Moment Matching
#JAGS model
model<-function(){
  Pc~dbeta(15,85)
  Pse~dbeta(3,9)
  OR.log~dnorm(-1.5,3)
  OR<-exp(OR.log)
  Pt<-exp(log(Pc/(1-Pc))+OR.log)/(exp(log(Pc/(1-Pc))+OR.log)+1)
  Qe.log~dnorm(0.6,6)
  Qe<-exp(Qe.log)/(exp(Qe.log)+1)
  
  #Future Sample
  Dc~dbin(Pc,200) 
  Dt~dbin(Pt,200)
  
  NB1<-Pc*(lambda*L*(1+Qe)/2-Ce)+
    (1-Pc)*lambda*L
  
  NB2<-Pse*Pt*(lambda*(L*(1+Qe)/2-Qse)-(Ct+Cse+Ce))+
    Pse*(1-Pt)*(lambda*(L-Qse)-(Ct+Cse))+
    (1-Pse)*Pt*(lambda*L*(1+Qe)/2-(Ct+Ce))+
    (1-Pse)*(1-Pt)*(lambda*L-Ct)
  INB<-NB2-NB1
}

Q<-30
var.post<-array()


#Calculate the posterior variance
for(i in 1:Q){
  #Sample the data
  Dc<-rbinom(1,200,quantile(Pc,i/(Q+1)))
  Dt<-rbinom(1,200,quantile(Pt,i/(Q+1)))
  #Running the JAGS model
  dataJags <- list(Dt=Dt,Dc=Dc,L=L,Qse=Qse,Ce=Ce,Ct=Ct,Cse=Cse,lambda=lambda)
  params<-c("INB","Pt","Pc")
  n.iter <- 10000
  n.burnin <- 1000
  post <- jags(dataJags, inits=NULL, params, model=model, n.iter=n.iter, 
                    n.burnin=n.burnin,DIC=FALSE)
  #Calculate posterior variance
  var.post[i]<-var(post$BUGSoutput$sims.list$INB)
}
#Calculating EVSI
fitted.full<-gam(INB~te(Pc,Pt))$fitted
sd.prepost<-sqrt(pr.var-mean(var.post))
#Perform moment matching
samp.prepost<-(fitted.full-mean(fitted.full))/sd(fitted.full)*
  sd.prepost+mean(fitted.full)

EVSI.MM<-mean(pmax(0,samp.prepost))-max(0,mean(samp.prepost))

###Full MCMC
###NOTE: This takes around 2.5 days to run.
mean.post<-array()
for(i in 1:1e+06){
  #Future data
  Dc<-rbinom(1,200,Pc[i])
  Dt<-rbinom(1,200,Pt[i])
  #Running the Bayesian model
  dataJags <- list(Dt=Dt,Dc=Dc,L=L,Qse=Qse,Ce=Ce,Ct=Ct,Cse=Cse,lambda=lambda)
  params<-c("INB","Pt","Pc")
  n.iter <- 10000
  n.burnin <- 1000
  vaccine <- jags(dataJags, inits=NULL, params, model=model,n.iter=n.iter,
                  n.burnin=n.burnin,DIC=FALSE)
  #Calculate the posterior mean
  mean.post[i]<-mean(vaccine$BUGSoutput$sims.list$INB)
}

EVSI.MCMC<-mean(pmax(0,mean.post))-max(0,mean(mean.post))
\end{verbatim}

\end{document}